\documentclass[12pt,reqno]{amsart}
\usepackage[english]{babel}
\usepackage{amsmath}
\usepackage[latin1]{inputenc}
\usepackage{amssymb}
\usepackage{amscd}
\usepackage{mathrsfs}
\usepackage[pdftex]{graphicx}
\usepackage{graphicx}

\newcommand{\N}[0]{\ensuremath{\mathbb{N}}}
\newcommand{\E}[0]{\ensuremath{\mathbb{E}}}
\newcommand{\F}[0]{\ensuremath{\mathbb{F}}}
\newcommand{\Pe}[0]{\ensuremath{\mathbb{P}}}
\newcommand{\R}[0]{\ensuremath{\mathbb{R}}}

\newtheorem{theo}{\sc{Theorem}}[section]
\newtheorem{lemm}{\sc{Lemma}}[section]

\newtheorem{cor}{\sc{Corollary}}[section]
\newtheorem{rmq}{\sc{Remark}}[section]

\title[]{Estimates for Solutions of a Low-Viscosity Kick-Forced Generalised Burgers Equation}
\author[]{Alexandre Boritchev}
\date{\today}
\begin{document}


\begin{abstract}
We consider a non-homogeneous generalised Burgers equation:
$$
\frac{\partial u}{\partial t} + f'(u)\frac{\partial u}{\partial x} - \nu \frac{\partial^2 u}{\partial x^2} = \eta^{\omega},\quad t \in \R,\ x \in S^1.
$$
Here, $\nu$ is small and positive, $f$ is strongly convex and satisfies a growth assumption, while $\eta^{\omega}$ is a space-smooth random "kicked" forcing term.
\\ \indent
For any solution $u$ of this equation, we consider the quasi-stationary regime, corresponding to $t \geq 2$. After taking the ensemble average, we obtain upper estimates as well as time-averaged lower estimates for a class of Sobolev norms of $u$. These estimates are of the form $C \nu^{-\beta}$ with the same values of $\beta$ for bounds from above and from below. They depend on $\eta$ and $f$, but do not depend on the time $t$ or the initial condition.
\end{abstract}

\maketitle{}

\section{Notation}
Consider a zero mean value smooth function $w$ on $S^1$.
For $p \in [1,+\infty]$, we denote its $L_p$ norm of by $\left|w\right|_p$. The $L_2$ norm will be denoted by $\left|w\right|$, and $\left\langle \cdot,\cdot\right\rangle$ stands for the $L_2$ scalar product. From now on, $L_p,\ p \in [1,+\infty]$ stands for the space of zero mean value functions in $L_p(S^1)$.
\\ \indent
For a nonnegative integer $n$ and $p \in [1,+\infty]$, $W^{n,p}$ stands for the Sobolev space of zero mean value functions $w$ on $S^1$ with the norm
\begin{equation} \nonumber
\left|w\right|_{n,p}=\left|w^{(n)}\right|_p,
\end{equation}
where
$$
w^{(n)}=\frac{d^n w}{dx^n}.
$$
In particular, $W^{0,p}=L_p$ for $p \in [1,+\infty]$. For $p=2$, we denote $W^{n,2}$ by $H^n$, and the corresponding norm is abbreviated as $\left\|w\right\|_n$.
\\ \indent
We recall a version of the classical Gagliardo-Nirenberg inequality (see \cite[p. 125]{Nir}).
\begin{lemm} \label{GN}
For a smooth zero mean value function $w$ on $S^1$,
$$
\left|w\right|_{\beta,r} \leq C \left|w\right|^{\theta}_{m,p} \left|w\right|^{1-\theta}_{q},
$$
where $m>\beta$, and $r$ is defined by
$$
\frac{1}{r}=\beta+\theta(\frac{1}{p} - m)+(1-\theta)\frac{1}{q},
$$
under the assumption that $\theta=\beta/m$ if $p=1$ or $p=+\infty$, and $\beta/m \leq \theta < 1$ otherwise. Here $C=C(m,p,q,\beta,\theta)>0$ is a constant.
\end{lemm}
\indent
For a smooth function $v(t,x)$ defined on $[0,+\infty) \times S^1$, $v_t$, $v_x$, and $v_{xx}$ mean respectively $\frac{\partial v}{\partial t}$, $\frac{\partial v}{\partial x}$, and $\frac{\partial^2 v}{\partial x^2}$.

\section{Introduction}
\indent
The generalised one-dimensional space-periodic Burgers equation
\begin{equation} \label{Burbegin}
\frac{\partial u}{\partial t} + f'(u) \frac{\partial u}{\partial x} - \nu \frac{\partial^2 u}{\partial x^2} = 0,\quad \nu>0
\end{equation}
(the classical Burgers equation corresponds to $f(u)=u^2$) appears in different domains of science, ranging from cosmology to traffic modelling (see \cite{BK}). It is sometimes called a viscous scalar conservation law. Historically, it has drawn most attention as a model for the Navier-Stokes equation (NSE).
Indeed, it has a nonlinear term analogous to the nonlinearity $(u \cdot \nabla) u$ in the incompressible NSE. The dissipation term in (\ref{Burbegin}) is also similar to the one in NSE. We note that the classical Burgers equation is explicitly solvable. This is done by the Cole-Hopf transformation (see \cite{BirCH}).
\\ \indent
In \cite{Bir}, A.Biryuk considered equation (\ref{Burbegin}) with $f$ strongly convex, i.e. satisfying
\begin{equation} \label{strconvex}
f''(x) \geq \sigma > 0,\quad x \in \R.
\end{equation}
He studied the behavior of the Sobolev norms of solutions $u$ for small values of $\nu$ and obtained the following estimates:
$$
\Vert u\Vert^2_m \leq C \nu^{-(2m-1)/2},\quad \frac{1}{T} \int_{0}^{T}{\Vert u\Vert^2_m} \geq c \nu^{-(2m-1)/2},\quad m \geq 1,\ \nu \leq \nu_0.
$$
Note that exponents of $\nu$ in lower and upper estimates are the same. The quantities $\nu_0$, $C$, $c$, and $T$ depend on the deterministic initial condition $u_0$ as well as on $m$. To get results independent from the initial data, a natural idea is to introduce random forcing and to estimate ensemble-averaged norms of solutions.
\\ \indent
In this article we consider (\ref{Burbegin}) with a random kick force in the right-hand side. In Section \ref{intro} we recall classical existence and uniqueness results and introduce the probabilistic setting needed to define the kick force. Then, we estimate from above the moments of the $W^{1,1}$ norm of $u$. These estimates, valid after a certain damping time, are proved using ideas similar to those in \cite{Kru}. Remarkably, this damping time and the estimate do not depend on the initial condition. This is the crucial result of this article.
\\ \indent
Next, in Sections \ref{lower} and \ref{upper}, this result allows us to obtain lower and upper estimates that are, up to taking the  ensemble average, of the same type as in \cite{Bir}, for time $t \geq 2$. These estimates will only depend on the function $f$ and the forcing. Let us emphasise that, for $t \geq 2$, we are in a  quasi-stationary regime: all estimates hold independently of the initial condition. In Section \ref{other}, we give some additional estimates for the Sobolev norms.
\\ \indent
In this paper, we use methods introduced by Kuksin in \cite{KukGAFA97,KukGAFA99}, and developed by Biryuk in \cite{Bir}.
\\ \indent
Equation (\ref{Burbegin}) with $\nu \ll 1$ is a popular one-dimensional model for the theory of hydrodynamic turbulence. In Section \ref{conclu}, we present an interpretation of our results in terms of this theory.

\section*{Acknowledgements}
\indent
First of all, I would like to thank my advisor S.Kuksin for formulation of the problem and guidance of my research. I would also like to thank A.Biryuk and K.Khanin for fruitful discussions. Finally, I am grateful to the faculty and staff at CMLS (Ecole Polytechnique) for their constant support during my PhD studies.

\section{Preliminaries} \label{intro}

\indent
In this section, we review properties of solutions of (\ref{Burbegin}) used in our proof.
\\ \indent
Physically, $t$ corresponds to the time variable, whereas $x$ corresponds to the one-dimensional space variable, and the constant $\nu>0$ to a viscosity coefficient. The real-valued function $u(t,x)$ is defined on $[0, +\infty) \times \R$ and is $L$-periodic in $x$. The function $f$ is $C^{\infty}$-smooth and strongly convex, i.e. it satisfies the condition (\ref{strconvex})
for some constant $\sigma$. Moreover, we assume that $f$, as well as its derivatives, has at most polynomial growth, i.e.
\begin{equation} \label{poly}
\forall m \geq 0,\ \exists n \geq 0,\ C_m>0:\quad |f^{(m)}(x)| \leq C_m (1+|x|)^n,\quad x \in \R,
\end{equation}
where $n=n(m)$. From now on, we fix $L=1$, which amounts to studying the problem on $[0, +\infty) \times S^1$.
We note that $L$-periodic solutions of (\ref{Burbegin}) with any $L$ reduce, by means of scaling in $x$, to $1$-periodic solutions with scaled $f$ and $\nu$.
\\ \indent
Since we are mostly interested in the asymptotics of solutions of (\ref{Burbegin}) as $\nu \rightarrow 0^+$, we assume that
\begin{equation} \nonumber
\nu \in (0,1].
\end{equation}
Moreover, it is enough to study the special case
\begin{equation} \label{mean0}
\int_{S^1}{u_0(y)dy}=0.
\end{equation}
Indeed, if the mean value of $u_0$ on $S^1$ equals $b$, we may consider
$$
v(t,x)=u(t,x+bt)-b.
$$
Then $v$ satisfies (\ref{mean0}) and is a solution of (\ref{Burbegin}) with $f(y)$ replaced with $g(y)=f(y+b)-by$.
\\ \indent
Given a $C^{\infty}$-smooth initial condition $u_0=u(0,\cdot)$, equation (\ref{Burbegin}) has a unique classical solution $u$,  $C^{\infty}$-smooth in both variables (see \cite[Chapter 5]{KrLo}). Condition (\ref{mean0}) implies that the mean value of a solution for (\ref{Burbegin}) vanishes identically in $t$.

Now provide each space $W^{n,p}(S^1)$ with the Borel $\sigma$-algebra. Consider a random variable $\zeta$ on a probability space $(\Omega, \F, \Pe)$ with values in $L^2(S^1)$, such that $\zeta^{\omega} \in C^{\infty}(S^1)$ for a.e. $\omega$. We suppose that $\zeta$ satisfies the following three properties.
\\ \indent
\textbf{(i) (Non-triviality)}
$$
\Pe(\zeta \equiv 0)<1.
$$
\indent
\textbf{(ii) (Finiteness of exponential moments for Sobolev norms)}
For every $m \geq 0$ there are constants $\alpha=\alpha(m)>0,\ \beta=\beta(m)$ such that
$$
\E \exp(\alpha \left\|\zeta\right\|^2_m) \leq \beta.
$$
In particular
$$
I_m=\E \left\|\zeta\right\|^{2}_m < +\infty,\quad \forall m \geq 0.
$$
\\ \indent
\textbf{(iii) (Vanishing of the expected value)}
$$
\E \zeta \equiv 0.
$$
\\ \indent
It is not difficult to construct explicitly
$\zeta$ satisfying \textbf{(i)-(iii)}. For instance
we could consider the real Fourier coefficients of $\zeta$, defined for $k>0$ by
\begin{equation} \label{Fourier}
a_k(\zeta)=\sqrt{2} \int_{S^1}{\cos(2 \pi kx) u(x)},\ b_k(\zeta)=\sqrt{2} \int_{S^1}{\sin(2 \pi kx) u(x)}
\end{equation}
as independent random variables with zero mean value and exponential moments tending to $1$ fast enough as $k \rightarrow + \infty$.
\\ \indent
Now let $\zeta_i$, $i \in \N$ be independent identically distributed random variables having the same distribution as $\zeta$. The sequence $(\zeta_i)_{i \geq 1}$ is a random variable, defined on a probability space which is a countable direct product of copies of  $\Omega$. From now on, this space will itself be called $\Omega$. The meaning of $\F$ and $\Pe$ changes accordingly.
\\ \indent
For $\omega \in \Omega$ and a time period $\theta>0$, the kick force $\eta^{\omega}$ is a $C^{\infty}$-smooth function in the variable $x$, with values in the space of distributions in the variable $t$, defined by
$$
\eta^{\omega}(x)=\sum_{i=1}^{+\infty}{\delta_{t=i \theta} \zeta_i^{\omega}(x)},
$$
where $\delta_{t=i \theta}$ denotes the Dirac measure at a time moment $i \theta$.
\\ \indent
The kick-forced version of (\ref{Burbegin}) corresponds to the case where, in the right-hand side, $0$ is replaced with the kick force.
This means that  for integers $i \geq 1$, at the moments $i\theta$ the solution $u(x)$ instantly increases by
the kick $\zeta_i^{\omega}(x)$, and that between these moments $u$ solves (\ref{Burbegin}). The equation is written as follows:
\begin{equation} \label{kickBurgers}
\frac{\partial u}{\partial t} + f'(u)\frac{\partial u}{\partial x} - \nu \frac{\partial^2 u}{\partial x^2} = \eta^{\omega}.
\end{equation}
Derivatives are taken in the sense of distributions.
\\ \indent
When studying solutions of (\ref{kickBurgers}), we will always assume that the initial condition $u_0=u(0,\cdot)$ is $C^{\infty}$-smooth.
Moreover, we normalise those solutions to be right-continuous in time at the kick moments $i\theta$. Such a solution is uniquely defined for a given value of $u_0$, for a.e. $\omega$.
\\ \indent
For a given initial condition $u_0$, the function $u(t,x)$ always will denote such a solution of (\ref{Burbegin}). The value of $u$ before the $i$-th kick will be denoted by $u(i \theta^-, \cdot)$, or shortly $u_i^-$. We will also use the notation $u_{i}=u(i \theta, \cdot)$ and denote the function $u(t,\cdot)$ by $u(t)$. Finally, for a solution of (\ref{kickBurgers}), we consider time derivatives at the kick moments in the sense of right-sided time derivatives. Those derivatives are right-continuous in time.
\\ \indent
Since space averages of the kicks vanish and $u_0(x)$ satisfies (\ref{mean0}), the space average of $u(t),\ t \geq 0$ vanishes identically. For the sake of simplicity, we normalise the kick period: from now on $\theta=1$.
\\ \indent
We observe that, since the kicks are independent and between the kicks (\ref{kickBurgers}) is deterministic, the solutions of (\ref{kickBurgers}) make a random Markov process. For details, see \cite{Kuk2}, where a kick force is introduced in a similar setting.

\bigskip
\indent
\textbf{Agreements.} All constants denoted $C$ with sub- or super-indexes are strictly positive. Unless otherwise stated, they depend only on $f$, on the distribution of the kicks, as well as on the parameters $a_1,\dots,a_k$ if they are denoted $C(a_1,\dots,a_k)$. $u$ always denotes a solution of (\ref{kickBurgers}) with any initial condition $u_0$. Averaging in ensemble corresponds to averaging in $\Pe$. All our estimates hold independently of the value of $u_0$.
\medskip \\ \indent
We observe that for every integer $i$ we have the following energy dissipation identity on the maximal kick-free intervals:
\begin{equation} \label{dissipation}
A_i = \left|u_i\right|^2 - \left|u_{i+1}^-\right|^2,
\end{equation}
where
\begin{equation} \label{Ai}
A_i = 2 \nu \int_{i}^{(i+1)}{\left\|u(t)\right\|_1^2\ dt}.
\end{equation}
Indeed, for any $t \in \left(i, i+1\right)$ $u$ satisfies
\begin{equation} \nonumber
2 \nu \left\|u(t)\right\|_1^2 = - 2 \nu \int_{S^1}{u u_{xx} dx} =  -2 \int_{S^1}{u f'(u) u_{x} dx} -2 \int_{S^1}{u u_{t} dx}.
\end{equation}
The first term on the right-hand side vanishes since its integrand is a full derivative.
The second term equals $-\frac{d}{dt}\left|u\right|^2$. Integrating in time we get (\ref{dissipation}).
We note that energy dissipation between kicks $A_i$ is always non-negative: energy can be added only at the kick points. We also note that an analogue of (\ref{dissipation}) holds on every kick-free time interval.
\\ \indent
\medskip
The following two lemmas are proved using the maximum principle in the same way as in \cite{Kru}.

\begin{lemm} \label{uxpos}
We have the estimate
$$
u_x(t,x) \leq 2 \sigma^{-1},\quad t \in [k+1/2,\ k+1),\ k \in \N,\ x \in S^1,
$$
where $\sigma$ is the constant in the assumption (\ref{strconvex}).
\end{lemm}

\textbf{Proof.}
Consider the equation (\ref{kickBurgers}) on the kick-free time interval $[0,1-\epsilon]$ for arbitrarily small $\epsilon$ and differentiate it once in space. We get
\begin{equation} \label{maxux}
\frac{\partial u_x}{\partial t} + f''(u) u^2_x + f'(u)\frac{\partial u_x}{\partial x} - \nu \frac{\partial^2 u_x}{\partial x^2} =0.
\end{equation}
Consider $v(t,x)=tu_x(t,x)$. For $t>0$, $v$ verifies
\begin{equation} \label{maxv}
\frac{\partial v}{\partial t} + t^{-1} (-v + f''(u) v^2) + f'(u)\frac{\partial v}{\partial x} - \nu \frac{\partial^2 v}{\partial x^2} =0.
\end{equation}
Now observe that, if $v>0$ somewhere on the domain $S_{\epsilon}=\left[0,1-\epsilon\right] \times S^1$, then $v$ attains its maximum $M$ on $S_{\epsilon}$ at a point $(t_1,x_1)$ such that $t_1>0$. At $(t_1,x_1)$ we have $\frac{\partial v}{\partial t} \geq 0$, $\frac{\partial v}{\partial x}=0$, and $\frac{\partial^2 v}{\partial x^2} \leq 0$.
Therefore, (\ref{maxv}) yields that
$$
t_1^{-1} [-v(t_1,x_1) + f''(u(t_1,x_1)) v^2(t_1,x_1)] \leq 0.
$$
Since, by (\ref{strconvex}), $f'' \geq \sigma > 0$, then
$$
-M+\sigma M^2 \leq 0,
$$
and therefore
$$
M \leq \sigma^{-1}.
$$
Thus we have proved that $v \leq \sigma^{-1}$ everywhere on $S_{\epsilon}$ for every $\epsilon>0$. In particular, by definition of $v$ and $S_{\epsilon}$, we get that
\begin{equation} \nonumber
u_x(t,x) \leq 2\sigma^{-1},\quad x \in S^1,\ t \in [1/2,1).
\end{equation}
Repeating the same argument on all the intervals $[k,k+1),\ k \in \N$ we get the lemma's assertion. $\qed$

\begin{lemm} \label{uxposbis}
There are constants $C',C$ such that
$$
\E \exp(C' \sup_{t \in [k,k+1)} \max u_x(t,\cdot)) \leq C,\quad k \geq 1.
$$
\end{lemm}

\textbf{Proof.}
Fix $k \geq 1$. Since the $W^{1,\infty}$ norm is dominated by the $H^2$ norm, then for $C'>0$ we get
\begin{align} \nonumber
\exp (C' u_x(k,x)) & \leq \exp (C' u_x(k^-,x) + C' \Vert\zeta_k\Vert_2),\quad x \in S^1.
\end{align}
The same inequality holds when we maximise in $x$. Now denote by $X_k$ the random variable
$$
\max u_x(k,\cdot).
$$
By Lemma \ref{uxpos} and Property \textbf{(ii)} of the kicks, for $C'=\alpha(2)$ we get
\begin{align}
\E \exp (C' X_k) & \leq \exp(2C' \sigma^{-1}) \E  \exp (C' \Vert\zeta_k\Vert_2) \leq C, \label{uxposbisaux}
\end{align}
for some constant $C$. Now consider the equation (\ref{maxux}). An application of the maximum principle to the function $u_x$, which cannot be negative everywhere, yields
\begin{align} \nonumber
\max u_x(t,\cdot) & \leq \max u_x(k,\cdot),\quad t \in \left[k,k+1\right).
\end{align}
Therefore, in (\ref{uxposbisaux}), we can replace $X_k$ by $\sup_{t \in [k,k+1)} \max u_x(t,\cdot)$. This proves the lemma's assertion. $\qed$


\begin{cor} \label{W11}
For the same $C',\ C$ as in Lemma \ref{uxposbis} we have
$$
\E \exp \Big( \frac{C'}{2} \sup_{t \in [k,k+1)} \left|u(t)\right|_{1,1} \Big) \leq C,\quad k \geq 1.
$$
\end{cor}

\textbf{Proof.}
Since the mean value of $u_x(t)$ is $0$, then
$$
\int_{S^1}{\left|u_x(t)\right|}=2 \int_{S^1}{\max(u_x(t),0)}.
$$
$\qed$

\begin{cor} \label{Lpupper}
For the same $C',\ C$ as in Lemma \ref{uxposbis} we have
$$
\E \exp(C' \sup_{t \in [k,k+1)} \left|u(t)\right|_{p}) \leq C,\quad k \geq 1,\ p \in [1,+\infty].
$$
Note that $C'$ and $C$ do not depend on $p$.
\end{cor}

\section{Lower estimates of $H^m$ norms} \label{lower}

For a solution $u$ of (\ref{kickBurgers}), the first quantity that we estimate from below is the expected value of
\begin{equation} \label{firstgoal}
\frac{1}{N} \int_{1}^{N+1}{\left\|u(t)\right\|_1^2}=\frac{1}{N} (2 \nu)^{-1}\sum_{i=1}^{N}{A_i},
\end{equation}
where $N$ is a fixed natural number chosen later, and $A_i$ is the same as in (\ref{Ai}).

\begin{lemm} \label{finitetime}
There exists a natural number $N \geq 1$, independent from $u_0$, such that
$$
\frac{1}{N} \int_{1}^{N+1}{\E \left\|u(s)\right\|_1^2} \geq C \nu^{-1}.
$$
\end{lemm}

\textbf{Proof.} For $N \geq 1$ we have
\begin{align} \nonumber
\E \left|u_{N+1}^-\right|^2 &\geq \E \Big( \left|u_{N+1}^-\right|^2-\left|u_1^-\right|^2 \Big)
\\ \nonumber
& = \E \sum_{i=1}^{N}{\Big( |u_{i+1}^-|^2 - |u_i|^2 \Big)}+ \E \sum_{i=1}^{N}{\Big( |u_{i}|^2 - |u_i^-|^2 \Big)}
\\ \nonumber
& = -\E \sum_{i=1}^{N}{A_i}+ \E \sum_{i=1}^{N}{\Big( |u_{i}^- + \zeta_i|^2 - |u_i^-|^2 \Big)}
\\ \nonumber
& = -\E \sum_{i=1}^{N}{A_i} + 2 \E \sum_{i=1}^{N}{\langle u_i^-,\zeta_i \rangle} + \E \sum_{i=1}^{N}{|\zeta_i|^2}.
\end{align} 
Since $\E \zeta_i \equiv 0$ (Property \textbf{(iii)} of the kicks), and $u_i^-$ and $\zeta_i$ are independent, then $\E \langle u_i^-,\zeta_i \rangle = 0$. Therefore, by (\ref{Ai}), we have
$$
\E \left|u_{N+1}^-\right|^2 \geq -2 \nu \E \int_{1}^{N+1}{\left\|u(s)\right\|_1^2} + 0 + NI_0.
$$
On the other hand, by Corollary \ref{Lpupper} ($p=2$) there is a constant $C_1$ such that
$$
\E \left|u_{N+1}^-\right|^2 \leq C_1.
$$
Consequently
$$
\frac{1}{N} \int_{1}^{N+1}{\E \left\|u(s)\right\|_1^2} \geq \frac{NI_0-C_1}{2N} \nu^{-1}.
$$
Choosing the smallest possible integer $N$ verifying
$$
N \geq \max \Big(1,\ \frac{C_1+1}{I_0} \Big),
$$
we get the lemma's assertion.\ $\qed$
\medskip \\ \indent
We have reached our first goal: estimating from below the expected value of (\ref{firstgoal}). Thus, we have a time-averaged lower estimate of the $H^1$ norm, which enables us to obtain similar estimates of $H^m$ norms for $m \geq 2$.

\begin{lemm} \label{finitetimebis}
We have
$$
\frac{1}{N} \int_{1}^{N+1}{\E \left\|u(s)\right\|_m^2} \geq C(m) \nu^{-(2m-1)},\quad m \geq 1, 
$$
where $N$ is the same as in Lemma \ref{finitetime}.
\end{lemm}

\textbf{Proof.}
This statement is already proved in the previous lemma for $m=1$, so we may assume that $m \geq 2$. By Lemma \ref{GN} and H{\"o}lder's inequality we have
\begin{equation} \label{Sobolev11}
\Big( \E \left\|u(s)\right\|^2_1 \Big)^{2m-1} \leq C'(m) \E \left\|u(s)\right\|_m^2 \Big(\E \left|u(s)\right|_{1,1}^2\Big)^{2m-2}.
\end{equation}
Since by Corollary \ref{W11}
$$
\E \left|u(s)\right|^2_{1,1} \leq K,\quad t \in [1,N+1],
$$
where $K>0$ is a constant, then, integrating (\ref{Sobolev11}) in time, we get
\begin{align} \nonumber
\frac{1}{N} \int_{1}^{N+1}{\E \left\|u(s)\right\|_m^2}  \geq & \frac{\int_{1}^{N+1} [\E (\left\|u(s)\right\|^2_1)]^{(2m-1)}}{N C'(m) K^{2m-2}}.
\end{align}
By H{\"o}lder's inequality,
\begin{equation} \nonumber
\int_{1}^{N+1} [\E (\left\|u(s)\right\|^2_1)]^{(2m-1)} \geq \Big(\int_{1}^{N+1}{\E \left\|u(s)\right\|^2_1} \Big)^{(2m-1)} N^{2-2m},
\end{equation}
and then
\begin{align} \nonumber
\frac{1}{N} \int_{1}^{N+1}{\E \left\|u(s)\right\|_m^2} \geq & \frac{\Big(\int_{1}^{N+1}{\E \left\|u(s)\right\|^2_1 \Big)^{(2m-1)}} N^{2-2m}}{N C'(m) K^{2m-2}}
\\ \nonumber
= & \frac{\Big( \frac{1}{N} \int_{1}^{N+1}{\E \left\|u(s)\right\|_1^{2}} \Big)^{(2m-1)}}{C'(m) K^{2m-2}}.
\end{align}
Now the assertion follows from Lemma \ref{finitetime}. $\qed$
\medskip \\ \indent
Since we impose no conditions on $u_0$, we can consider a different positive integer "starting time". We may also consider a different averaging time interval of length $T \geq N$. Finally, we obtain a general result for a non-integer starting time $t \geq 1$ by considering the maximal interval $[m_1,m_2] \subset [t,t+T]$ such that $m_1$ and $m_2$ are positive integers.

\begin{theo} \label{finalexp}
We have
$$
\frac{1}{T} \int_{t}^{t+T}{\E \left\|u(s)\right\|_m^2} \geq \frac{C(m)}{4} \nu^{-(2m-1)},\quad t \geq 1,\ T \geq N+1,\ m \geq 1,
$$
where $N$ and $C(m)$ are the same as in Lemma \ref{finitetimebis}.
\end{theo}

\section{Upper estimates of $H^m$ norms} \label{upper}
\indent
To estimate from above a Sobolev norm $\left\|u\right\|_m,\ m \geq 1,$ of a solution $u$ for (\ref{kickBurgers}), we differentiate between the kicks the quantity $\left\|u(t)\right\|_m^2$.
\\ \indent
Denote by $B(u)$ the nonlinearity $2 f'(u)u_x$, and by $L$ the operator $-\partial_{xx}$.
Integrating by parts, we get
\begin{align} \nonumber
\frac{d}{dt} \left\|u\right\|_m^2 & = 2 \left\langle u^{(m)}, u_t^{(m)} \right\rangle
\\ \label{diffum2}
& = -2 \nu \left\|u\right\|_{m+1}^2  - \left\langle L^{m}u, B(u) \right\rangle.
\end{align}
\indent
We will need a standard estimate for the nonlinearity $\left\langle L^m u, B(u)\right\rangle$.

%

\begin{lemm} \label{lmubuinfty}
For a zero mean value smooth function $w$ such that $\left|w\right|_{\infty} \leq M$, we have
$$
\left| \left\langle L^m w, B(w) \right\rangle\right| \leq C \left\|w\right\|_m \left\|w\right\|_{m+1},\quad m \geq 1,
$$
with $C$ satisfying
\begin{equation} \label{polyestimate}
C \leq C_m (1+M)^n,
\end{equation}
where $C_m$, as well as the natural number $n=n(m)$, depend only on $m$.
\end{lemm}

\textbf{Proof.}
Let $C'$ denote various positive constants satisfying an estimate of the type (\ref{polyestimate}). Then we have
\begin{align} \nonumber
\left|\left\langle L^m w, B(w)\right\rangle\right| = & 2 \left|\left\langle w^{(2m)}, (f(w))^{(1)} \right\rangle\right|
\\ \nonumber
= & 2 \left|\left\langle w^{(m+1)}, (f(w))^{(m)} \right\rangle\right|
\\ \nonumber
\leq & C' \sum_{k=1}^m\ \sum_{\substack{1 \leq a_1 \leq \dots \leq a_k \leq m \\ a_1+ \dots+a_k = m}} \int_{S^1}{\left| w^{(m+1)} w^{(a_1)} \dots w^{(a_k)} f^{(k)}(w) \right|}
\\ \nonumber
\leq & C' \left|f\right|_{C^m[-M,M]} \sum_{k=1}^m\ \sum_{\substack{1 \leq a_1 \leq \dots \leq a_k \leq m \\ a_1+ \dots+a_k = m}} \int_{S^1} | w^{(a_1)} \dots 
\\ \nonumber
& \dots w^{(a_k)} w^{(m+1)} |.
\end{align}
By (\ref{poly}), $\left|f\right|_{C^m[-M,M]}$ satisfies an estimate of the type (\ref{polyestimate}).
By H{\"o}lder's inequality, we obtain that
\begin{align} \nonumber
\left|\left\langle L^m w, B(w) \right\rangle\right| \leq & C' \left\|w\right\|_{m+1} \sum_{\substack{1 \leq a_1 \leq \dots \leq a_k \leq m \\ a_1+ \dots+a_k = m}} \Big( \left|w^{(a_1)}\right|_{2m/a_1} \dots 
\\ \nonumber
& \dots \left|w^{(a_k)}\right|_{2m/a_k} \Big).
\end{align}
Finally, the Gagliardo-Nirenberg inequality yields
\begin{align} \nonumber
\left|\left\langle L^m w, B(w)\right\rangle\right| \leq & C' \left\|w\right\|_{m+1} \sum_{k=1}^m\ \sum_{\substack{1 \leq a_1 \leq \dots \leq a_k \leq m \\ a_1+ \dots+a_k = m}}
\\ \nonumber
&  \Big[ (\left\|w\right\|_m^{a_1/m} \left|w\right|_{\infty}^{(m-a_1)/m} ) \dots (\left\|w\right\|_m^{a_k/m} \left|w\right|_{\infty}^{(m-a_k)/m}) \Big]
\\ \nonumber
\leq & C' \left|w\right|_{\infty}^{m-1} \left\|w\right\|_m \left\|w\right\|_{m+1}
\\ \nonumber
\leq & C' \left\|w\right\|_m \left\|w\right\|_{m+1},
\end{align}
which proves the lemma's assertion.\ $\qed$

\begin{theo} \label{upperkickm}
For any natural numbers $m,n$ we have
$$
\E ( \sup_{t \in [k,k+1)} \left\|u(t)\right\|^{n}_m) \leq C(m,n)  \nu^{-(2m-1)n/2},\quad k \geq 2.
$$
\end{theo}

\textbf{Proof.} 
Fix $k \geq 2$ and $m \geq 1$.
In this proof, $\Theta$ denotes various positive random constants which depend on $m$, such that all their moments are finite, and $C$ denotes various positive deterministic constants, depending only on $m$.
\\ \indent
We begin by noting that Corollary \ref{W11} and Property \textbf{(ii)} of the kicks imply the inequalities
\begin{equation} \label{uxposinside}
|u(t)|_{1,1}, \left\| \zeta_k \right\|_m \leq \Theta,\quad t \in [k-1,k+1).
\end{equation}
\\ \indent
We claim that when $\left\|u\right\|^2_m$ is too large, it decreases at least as fast as a solution of the differential equation
$$
y'+(2m-1) y^{2m/(2m-1)}=0,
$$
i.e. as $t^{-(2m-1)}$. More precisely, we want to prove that for
\\
$t \in [k-1,k+1)$ we have
\begin{align} \nonumber
&\left\|u(t)\right\|^2_m \geq \Theta_1 \nu^{-(2m-1)} \Longrightarrow
\\ \label{decrm}
&\ \frac{d}{dt}\left\|u(t)\right\|^2_m \leq - (2m-1) \left\|u(t)\right\|^{4m/(2m-1)}_m,
\end{align}
where $\Theta_1$ is a random positive constant, chosen later. Random constants $\Theta$ below do not depend on $\Theta_1$.
\\ \indent
Indeed, assume that
\begin{equation} \label{C1assum}
\left\|u(t)\right\|^2_m \geq \Theta_1 \nu^{-(2m-1)}.
\end{equation}
We begin by observing that by Lemma \ref{GN} we have
\begin{equation} \nonumber
\left\|u\right\|_m \leq C \left\|u\right\|_{m+1}^{(2m-1)/(2m+1)} \left|u\right|_{1,1}^{2/(2m+1)},
\end{equation}
and hence
\begin{align} \nonumber
\left\|u\right\|_{m+1} & \geq C \left|u\right|_{1,1}^{-2/(2m-1)} \left\|u\right\|_m^{(2m+1)/(2m-1)}
\\  \label{um+12}
& \geq \Theta^{-1} \left\|u\right\|_m^{(2m+1)/(2m-1)}
\end{align}
(we used (\ref{uxposinside})). Now, (\ref{diffum2}), (\ref{uxposinside}), and Lemma \ref{lmubuinfty} imply that
\begin{align} \nonumber
\frac{d}{dt} & \left\|u\right\|_m^2 \leq  - 2 \nu \left\|u\right\|_{m+1}^2 + \Theta \left\|u\right\|_m \left\|u\right\|_{m+1}
\\ \label{dispm}
 = &(-2 \nu \left\|u\right\|_{m+1}^{2/(2m+1)} + \Theta \left\|u\right\|_m \left\|u\right\|_{m+1}^{-(2m-1)/(2m+1)}) \left\|u\right\|_{m+1}^{4m/(2m+1)}.
\end{align}
Combining (\ref{dispm}) and (\ref{um+12}), we get
\begin{align} \nonumber
\frac{d}{dt} \left\|u\right\|_m^2 \leq & (-2 \nu \left\|u\right\|_{m+1}^{2/(2m+1)} + \Theta) \left\|u\right\|_{m+1}^{4m/(2m+1)}.
\end{align}
Therefore, by (\ref{um+12}) and (\ref{C1assum}) we have
\begin{align} \nonumber
\frac{d}{dt} \left\|u\right\|_m^2 \leq & \Big(-\nu  \Theta^{-1} \left\|u\right\|_{m}^{2/(2m-1)} 
+ \Theta \Big) \left\|u\right\|_{m+1}^{4m/(2m+1)}
\\ \nonumber
\leq & \left(- \Theta^{-1} \Theta_1^{1/(2m-1)}+ \Theta \right) \left\|u\right\|_{m+1}^{4m/(2m+1)}.
\end{align}
Now we choose $\Theta_1$ in such a way that the quantity in the parentheses is negative. Under this assumption, we get from (\ref{um+12}) that
\begin{align} \nonumber
\frac{d}{dt} \left\|u\right\|_m^2 \leq & \left(-\Theta^{-1} \Theta_1^{1/(2m-1)}+ \Theta \right) \Theta^{-1} \left\|u\right\|_{m}^{4m/(2m-1)}.
\end{align}
This relation implies (\ref{decrm}) if we choose for $\Theta_1$ a sufficiently big random constant with all moments finite.
\\ \indent
Now we claim that
\begin{equation} \label{decrmcor}
\left\|u^{-}_{k}\right\|^2_m \leq \Theta_2 \nu^{-(2m-1)},
\end{equation}
where
$$
\Theta_2=\max(\Theta_1,\ 1)
$$
has finite moments.
Indeed, if $\left\|u(t)\right\|^2_m \leq \Theta_1 \nu^{-(2m-1)}$ for some $t \in \left[k-1,k\right)$, then (\ref{decrm}) ensures that $\left\|u(t)\right\|^2_m$ remains under this threshold up to $t=k^-$. Otherwise, we consider the function
$$
y(t)=\left\|u(t)\right\|^{-2/(2m-1)}_m,\quad t \in \left[k-1,k\right).
$$
By (\ref{decrm}), since $\left\|u(t)\right\|^2_m > \Theta_1 \nu^{-(2m-1)}$, $y(t)$ increases at least as fast as $t$. Indeed,
\begin{align} \nonumber
\frac{d}{dt} y(t) & = -\frac{1}{2m-1}\  \Big( \left\|u(t)\right\|^2_m \Big)^{-2m/(2m-1)} \frac{d}{dt} \left\|u(t)\right\|^2_m
\\ \nonumber
& \geq \frac{1}{2m-1} \left\|u(t)\right\|^{-4m/(2m-1)}_m (2m-1)  \left\|u(t)\right\|^{4m/(2m-1)}_m 
\\ \nonumber
& \geq 1. 
\end{align}
Therefore $\left\|y(k^-)\right\|^2_m \geq 1$. Since $\nu \leq 1$, then in this case we also have (\ref{decrmcor}).
\\ \indent
In exactly the same way, using (\ref{uxposinside}), we obtain that for $t \in [k,k+1)$,
\begin{align} \nonumber
\left\|u(t)\right\|^2_m & \leq \max(\Theta_2\nu^{-(2m-1)},\left\|u(k)\right\|_m^2)
\\ \nonumber
& \leq \max \Big[\Theta_2, \Big(\Theta+\sqrt{\Theta_2}\Big)^2 \Big] \nu^{-(2m-1)}
\\ \nonumber
& \leq \Big(\Theta+\sqrt{\Theta_2}\Big)^2 \nu^{-(2m-1)}.
\end{align}
Therefore $\left\|u(t)\right\|^2_m \nu^{2m-1}$ is uniformly bounded by $\Big(\Theta+\sqrt{\Theta_2}\Big)^2$ for $t \in [k,k+1)$. Since all moments of this random variable are finite,  the lemma's assertion is proved. $\qed$

\section{Estimates of other Sobolev norms.} \label{other}
\indent
The results in the three previous sections enable us to find upper and lower estimates for a large class of Sobolev norms. 
Unfortunately, while lower estimates extend to the whole Sobolev scale for $m \geq 0$ and $p \in [1,+\infty]$, there is a gap, corresponding to the case $m \geq 2$ and $p=1$, for upper estimates.

\begin{lemm} \label{upperkickwmp}
For $m \in \left\lbrace 0,1 \right\rbrace $ and $p \in [1,+\infty]$, or for $m \geq 2$ and $p \in (1,+\infty]$, we have
$$
\Big( \E \sup_{t \in [k,k+1)} \left|u(t)\right|^{n}_{m,p} \Big)^{1/n} \leq C(m,p,n)  \nu^{-\gamma},\quad n \geq 1,\ k \geq 2.
$$
Here and later on,
$$
\gamma=\gamma(m,p)=\max \Big( 0,\ m-\frac{1}{p} \Big).
$$
\end{lemm}

\textbf{Proof.}
We begin by considering the case $m=1$ and $p \in [2,+\infty]$. Since by Lemma \ref{GN} we have
$$
\left|u(t)\right|_{m,p} \leq C(m,p) \left\|u(t)\right\|^{1-\theta}_{m} \left\|u(t)\right\|^{\theta}_{m+1},
$$
where
$$
\theta=\frac{1}{2}-\frac{1}{p},
$$
then Theorem \ref{upperkickm} and H{\"o}lder's inequality yield the wanted result.
\\ \indent
The case $m=1$ and $p \in [1,2)$ is proved in exactly the same way, by combining Corollary \ref{W11} and Theorem \ref{upperkickm} ($m=1$). The same method is used to prove the case $m \geq 2$ and $p \in (1,2)$, combining the case $p \in [2,+\infty]$ for a big enough value of $m$ and Corollary \ref{W11}. Unfortunately, it cannot be applied for $m \geq 2$ and $p=1$, because Lemma \ref{GN} only allows us to estimate a $W^{n,1}$ norm from above by other $W^{n,1}$ norms.
\\ \indent
Finally, the case $m=0$ follows from Corollary \ref{Lpupper}.
$\qed$
\\ \smallskip \\ \indent
The first norm that we estimate from below is the $L_2$ norm.

\begin{lemm} \label{L2lower}
We have
\begin{align} \nonumber
\Big(\int_{k}^{k+1} \E |u(s)|^{2} \Big)^{1/2} \geq C,\quad k \geq 2.
\end{align}
\end{lemm}

\textbf{Proof.}
Using Properties \textbf{(i)} and \textbf{(iii)} of the kicks ($u_k^-$ and $\zeta_k$ being independent), we get
\begin{align} \nonumber
\E \left|u_k^+\right|^2 &= \E \left|u_k^-\right|^2+2 \E \left\langle u_k^-,\ \zeta_k\right\rangle + \E \left|\zeta_k\right|^2
\\ \nonumber
&= \E \left|u_k^-\right|^2+\E \left|\zeta_k\right|^2 \geq I_0.
\end{align}
On the other hand, by Theorem \ref{upperkickm} we have
$$
\E \left\|u(t)\right\|_1^2 \leq C' \nu^{-1},\quad t \in (k,k+1).
$$
Since
$$
\frac{d}{dt} \left|u(t)\right|^2 = -2 \nu \left\|u(t)\right\|_1^2,\quad t \in (k,k+1),
$$
then, integrating in time and setting
$$
d=\min \Big( 1,\frac{I_0}{4C'} \Big),
$$
we obtain that, for $s \in [k,k+d]$,
$$
\E |u(s)|^{2} \geq \E \left|u_k^+\right|^2 - 2(s-k) C' \geq I_0-2C' d \geq \frac{I_0}{2}.
$$
Therefore
\begin{align} \nonumber
\int_{k}^{k+1} \E |u(s)|^{2} & \geq \min \Big( \frac{I_0}{2},\ \frac{I_0^2}{8C'} \Big)>0,
\end{align}
which proves the lemma's assertion. $\qed$
\bigskip
\\ \indent
Now we can study the case $m=0$ and $p \in [1,+\infty]$.

\begin{cor} \label{Lplower}
We have
\begin{align} \nonumber
\Big( \int_{k}^{k+1} \E|u(s)|_{p}^{2} \Big)^{1/2} \geq C,\quad k \geq 2,\ p \in [1,+\infty],
\end{align}
where $C$ does not depend on $p$.
\end{cor}

\textbf{Proof.}
It suffices to prove the inequality for $p=1$.
Using H{\"o}lder's inequality and integrating in time and in ensemble, and then using the Cauchy-Schwarz inequality, we get
\begin{align} \nonumber
\int_{k}^{k+1} \E |u|_1^{2} & \geq \int_{k}^{k+1} \E |u|^{4} |u|_{\infty}^{-2}
\\ \nonumber
& \geq \Big( \int_{k}^{k+1} \E |u|^{2} \Big)^2 \Big(  \int_{k}^{k+1} \E |u|_{\infty}^{2} \Big)^{-1}.
\end{align}
Lemma \ref{L2lower} and Corollary \ref{Lpupper} ($p=+\infty$) complete the proof. $\qed$
\bigskip
\\ \indent
Since  the $W^{1,1}$ norm dominates the $L_{\infty}$ norm, we get

\begin{cor} \label{W11lower}
We have
\begin{align} \nonumber
\Big( \int_{k}^{k+1} \E |u(s)|^2_{1,1}(t) \Big)^{1/2} \geq C,\quad k \geq 2.
\end{align}
\end{cor}
The cases $m \geq 2$ and $m=1,\ p \geq 2$ follow from Lemma \ref{finitetime}  and Lemma \ref{GN} by interpolation in the same way as Lemma \ref{finitetimebis}, for $p>1$. The case $p=+\infty$ follows from the case $p=1$, since $|u|_{m,1} \geq |u|_{m-1,\infty}$, and $\gamma(m,1)=\gamma(m-1,+\infty)$.

\begin{lemm} \label{lowerbig}
If either $m \geq 2$ and $p \in [1,+\infty]$, or $m=1$ and $p \in [2,+\infty]$, then
$$
\Big(\frac{1}{T} \int_{t}^{t+T}{\E \left|u(s)\right|_{m,p}^2} \Big)^{1/2} \geq C(m,p) \nu^{-\gamma},\quad t \geq 1,\ T \geq N+1,
$$
where $N$ is the same as in Lemma \ref{finitetime}.
\end{lemm}

Now it remains to deal with the case $m=1$ and $p \in (1,2)$. 

\begin{lemm} \label{lowerm1p12}
For $p \in (1,2)$ we have
$$
\Big(\frac{1}{T} \int_{t}^{t+T}{\E \left|u(s)\right|_{1,p}^2} \Big)^{1/2} \geq C(p) \nu^{-\gamma},\quad t \geq 2,\ T \geq N+1,
$$
where $N$ is the same as in Lemma \ref{finitetime}. Note that here, $\gamma=1-1/p$.
\end{lemm}

\textbf{Proof.}
In the proof of this lemma, $C'(p)$ denotes various positive constants depending only on $p$. By H{\"o}lder's inequality in space we have
$$
\left\|u(s)\right\|^2_{1} \leq \left|u(s)\right|^{p}_{1,p} \left|u(s)\right|^{(2-p)}_{1,\infty}.
$$
Therefore, using H{\"o}lder's inequality in time and in ensemble, as well as Lemma \ref{upperkickwmp}, we get
\begin{align} \nonumber
\frac{1}{T} \int_{t}^{t+T}{\E \left\|u(s)\right\|_{1}^2} \leq & \Big(\frac{1}{T} \int_{t}^{t+T} \E \left|u(s)\right|^{2}_{1,\infty} \Big)^{(2-p)/2} \cdot
\\ \nonumber
& \Big(\frac{1}{T} \int_{t}^{t+T} \E \left|u(s)\right|^{2}_{1,p} \Big)^{p/2} 
\\ \nonumber
\leq &\ C'(p) \nu^{(p-2)} \Big(\frac{1}{T} \int_{t}^{t+T} \E \left|u(s)\right|^{2}_{1,p} \Big)^{p/2}.
\end{align}
Furthermore, Lemma \ref{finitetime} implies that
\begin{align} \nonumber
\frac{1}{T} \int_{t}^{t+T} \E \left|u(s)\right|^{2}_{1,p} & \geq C'(p) \Big( \nu^{(2-p)} \frac{1}{T} \int_{t}^{t+T}{\E \left\|u(s)\right\|_{1}^2} \Big)^{2/p}
\\ \nonumber
& \geq C'(p) \Big( \nu^{(2-p)} \nu^{-1} \Big)^{2/p}
\\ \nonumber
& \geq C'(p) \nu^{{-(2p-2)/p}}.
\end{align}
$\qed$

\begin{rmq}
Upper estimates for
$$
\Big(\frac{1}{T} \int_{t}^{t+T} \E |u(s)|^{n}_{m,p} \Big)^{1/n},\quad n \geq 2
$$
follow from the lemmas above and H{\"o}lder's inequality.
\end{rmq}

\section{Conclusion} \label{conclu}

Putting together the estimates that we have obtained, we formulate our main result.

\begin{theo} \label{general}
For $m \in \left\lbrace 0,1 \right\rbrace $ and $p \in [1,+\infty]$, or for $m \geq 2$ and $p \in (1,+\infty]$, we have
\begin{align} \label{uppergeneral}
\Big( \E \sup_{t \in [k,k+1)} |u(t)|^{n}_{m,p} \Big)^{1/n} \leq C(m,p,n) \nu^{-\gamma},\quad n \geq 1,\ k \geq 2.
\end{align}
Moreover, there is an integer $N' \geq 1$ such that, for $m \geq 0$ and $p \in [1,+\infty]$, we have
\begin{align} \label{lowergeneral}
\Big(\frac{1}{T} \int_{t}^{t+T} \E |u(s)|^{n}_{m,p} \Big)^{1/n} \geq C(m,p) \nu^{-\gamma},\ \ n \geq 2,\ t \geq 2,\ T \geq N'.
\end{align}
In both inequalities
$$
\gamma=\max \Big( 0,\ m-\frac{1}{p}\ \Big).
$$
\end{theo}
\smallskip
For a solution $u$ of (\ref{kickBurgers}), we have obtained asymptotic estimates for expectations of a large class of Sobolev norms. The power of $\nu$ is clearly optimal except for $m \geq 2$ and $p=1$, since it coincides for upper and lower estimates: we are in a \textit{quasi-stationary regime}.  Let us stress again that the upper bound $t = 2$ for the time needed for a quasi-stationary regime to be established has no dependence on $u_0$. The condition $t \geq T_0$ for some time $T_0 \geq 1$ is necessary: we need damping if $u_0$ is large and injection of energy at a kick point if $u_0$ is small.
\\ \indent
Now put $\hat{u}^k=a_k(u)+i b_k(u)$ (see (\ref{Fourier})). For $t \geq 2$ and $T$ big enough (see Theorem \ref{general}), consider the averaged quantities
$$
F_{s,\theta}=\frac{1}{T} \int_t^{t+T} \frac{\sum_{k \in I(s,\theta)}{\E |\hat{u}^k|^2(\tau)}}{\sum_{k \in I(s,\theta)}{1}},\ \ \ \ s,\theta>0,
$$
where $I(s,\theta)=[\nu^{-s+\theta},\nu^{-s-\theta})$. In the same way as in \cite[formulas (1.6)-(1.8)]{Bir}, the inequalities (\ref{uppergeneral}-\ref{lowergeneral}) yield
\begin{align} \label{spectr1}
&F_{s,\theta} \leq C \nu^{2s}
\\ \label{spectr2}
&F_{s,\theta} \leq C(m) \nu^{2+2m(s-1-\theta)},\ \ m>0,\ s>1+\theta
\\ \label{spectr3} 
&F_{1,\theta}>C \nu^{2+2 \theta}
\end{align}
for $\nu \leq \nu(\theta)$ with some $\nu(\theta)>0$. These results have some consequences for the energy spectrum of $u$.
\\ \indent
Indeed, relation (\ref{spectr2}) implies that the energy of the $k$-th Fourier mode, $E_k=\frac{1}{2T} \int_t^{t+T} \E |\hat{u}^k|^2$, averaged around $k=l$, where $l \gg \nu^{-1}$, decays faster than any negative degree of $l$. On the other hand, by (\ref{spectr1}) and  (\ref{spectr3}), the energy $E_k$, averaged around $k=\nu^{-1}$, behaves as $k^{-2}$. That is, the interval $k \in (\nu^{-1}, +\infty)$ is the \textit{dissipation range}, where the energy $E_k$ decays fast.
\\ \indent
As the force $\eta$ is smooth in $x$, then the energy is injected at frequencies $k \sim 1$. The estimate (\ref{spectr1}) readily implies that the energy $E=\sum{E_k}$ of a solution $u$ is supported, when $\nu \rightarrow 0$, by any interval $(0,\nu^{-\gamma})$, $\gamma>0$. That is, the \textit{energy range} of the solution $u$ is the interval $(0,\nu^0]$ (see \cite{Fri}). 
\\ \indent
The complement to the energy and dissipation ranges is the \textit{inertial  range} $(\nu^{0},\nu^{-1})$. At $k \sim \nu^{-1}$ we have $E_k \sim k^{-2}$. It is plausible that in this range $E_k$ decays algebraically; possibly $E_k \sim k^{-2}$. The study of the energy spectrum of solutions $u$ in the inertial range is one of the objectives of our future research.
\\ \indent
We recall that the behavior of the energy spectrum $E_k$ of turbulent fluid of the form "some negative degree of $k$ in the inertial range, followed by fast decay in the dissipation range" is suggested by the Kolmogorov theory of turbulence (see \cite{Fri}). Our results (following those of A.Biryuk in \cite{Bir}) show that for the "burgulence" (described by the Burgers equation, see \cite{BK}) the dissipation range is $(\nu^{-1},+\infty)$ and suggest that the power-law in the inertial range is $E_k \sim k^{-2}$.
\\ \indent
We also see that for $\nu \rightarrow 0^+$, solutions $u$ display intermittency-type behavior (see \cite[Chapter 8]{Fri}). Indeed, in the quasi-stationary regime, up to averaging in time and in ensemble, $\max_{x \in S^1} u_x
 \sim 1$, whereas $\int_{S^1}{u_x^2} \sim \nu^{-1}$. Thus, typically $u$ has large negative gradients on a small subset of $S^1$, and small positive gradients on a large subset of $S^1$.
\\ \indent
In a future paper, we will look at the same problem with the kick force replaced by a spatially smooth white noise in time (see \cite{EKMS} for a possible definition). This problem is, heuristically, the limit case of the kick-forced problem with more and more frequent appropriately scaled kicks.

\bigskip
\begin{center}
Alexandre Boritchev
\\
Centre de Math{\'e}matiques Laurent Schwartz
\\
Ecole Polytechnique, Route de Saclay
\\
91128 Palaiseau Cedex,  France
\\
E-mail: boritchev@math.polytechnique.fr
\end{center}

\end{document}